\documentclass[letterpaper, 10 pt, conference]{ieeeconf}  

\IEEEoverridecommandlockouts                              
\usepackage{graphics} 
\usepackage{epsfig} 
\usepackage{amsmath} 
\usepackage{amssymb}  
\usepackage[dvipsnames]{xcolor}

\usepackage{stmaryrd} 
\usepackage{bbm} 
\usepackage{algorithmic}
\usepackage[linesnumbered,ruled,vlined]{algorithm2e}
\usepackage{url}
\usepackage{caption}
\usepackage{subcaption}
\usepackage{booktabs}
\usepackage{dsfont}
\usepackage{mathtools}

\newcommand{\bsgamma}{{\boldsymbol{\gamma}}}
\newcommand{\bsX}{\boldsymbol{X}}
\newcommand{\bsY}{\boldsymbol{Y}}
\newcommand{\bsZ}{\boldsymbol{Z}}
\newcommand{\bsB}{\boldsymbol{B}}
\newcommand{\bsA}{\boldsymbol{A}}
\newcommand{\bsalpha}{\boldsymbol{\alpha}}
\newcommand{\RR}{\mathbb{R}}
\newcommand{\EE}{\mathbb{E}}
\newcommand{\bsG}{\boldsymbol{\rm G}}
\newcommand{\diag}{\mathrm{diag}}

\DeclareMathOperator*{\argmin}{arg\,min}

\newtheorem{theorem}{Theorem}[section]
\newtheorem{remark}[theorem]{Remark}

\newtheorem{definition}[theorem]{Definition}

\newtheorem{proposition}[theorem]{Proposition}

\usepackage[colorlinks]{hyperref}

\definecolor{burgundy}{rgb}{0.5,0.0, 0.13}
\hypersetup{colorlinks,
citecolor=MidnightBlue,
citebordercolor=burgundy,
linkcolor=burgundy,
urlcolor=MidnightBlue
}

\title{\LARGE \bf
Multi-population Mean Field Games with Multiple Major Players:\\ Application to Carbon Emission Regulations
}

\author{G\"ok{\c c}e Dayan{\i}kl{\i} and Mathieu Lauri{\`e}re 
\thanks{Department of Statistics,
  University of Illinois at Urbana-Champaign, 
  Champaign, IL 61820, USA 
        {\tt\small gokced@illinos.edu}}%
\thanks{NYU-ECNU Institute of Mathematical Sciences at NYU Shanghai; Shanghai Frontiers Science Center of Artificial Intelligence and Deep Learning; NYU Shanghai, 567 West Yangsi Road, Shanghai, 200126, People’s Republic of China
        {\tt\small mathieu.lauriere@nyu.edu}}%
}

\begin{document}

\maketitle
\thispagestyle{empty}
\pagestyle{empty}

\begin{abstract}

 In this paper, we propose and study a mean field game model with multiple populations of minor players and multiple major players, motivated by applications to the regulation of carbon emissions. Each population of minor players represent a large group of electricity producers and each major player represents a regulator. We first characterize the minor players equilibrium controls using forward-backward differential equations, and show existence and uniqueness of the minor players equilibrium. We then express the major players' equilibrium controls through analytical formulas given the other players' controls. Finally, we then provide a method to solve the Nash equilibrium between all the players, and we illustrate numerically the sensitivity of the model to its parameters.

\end{abstract}

\section{Introduction} 
\label{sec:introduction}

Mean field games (MFGs) have been introduced to study large populations of agents that make decisions in a strategic manner, while interacting symmetrically~\cite{lasry2007mean,huang2006large}. They provide a simple and mathematically rigorous framework to approximate Nash equilibria in games with finite but large number of players. So far, most of the literature has focused on the case of a single population of homogeneous agents~\cite{cardaliaguet2018mean,bensoussan2013mean,carmona2018probabilistic}. Several extensions have been studied, in particular multi-population MFGs, in which there are several types of players with possibly different dynamics and cost functions~\cite{cirant2015multi,bensoussan2018meanmulti}, major-minor MFGs, in which a few  players have a macroscopic influence~\cite{nourian2013epsilon,carmonazhu2016probabilistic,advertisement}, and Stackelberg MFGs, in which one player chooses her control before a mean field population reacts to it~\cite{elie2019tale,dayanikli2023machine}. In the literature, many potential applications of MFGs on stylized models have been proposed, such as price formation~\cite{lachapelle2016efficiency,gomessaude2021meanprice}, systemic risk~\cite{carmonafouquesun2015mean}, epidemic management~\cite{elie2020contact,aurell2022optimal}, energy management~\cite{alasseur2020extended,gueant_2010,sircar_2017}. 

Climate change and its consequences are one of the main challenges that humankind currently faces. An important cause is carbon emissions, which are largely due to the electricity production. When there is a large number of producers or consumers, a mean field approximation can help to model the carbon emission and regulation in an efficient way. For example, \cite{malhame_2017} proposes an MFG to model climate change negotiations between countries that interact through a carbon emission market, \cite{carbon_StackelbergMFG} analyzes the problem of carbon tax implementation by using a Stackelberg single population mean field game, and~\cite{shrivats2021principal} studies the implementation of renewable energy certificates. \cite{escribe2023mean} proposes an MFG model for investments in renewable with common noise.

In this work, we are interested in modeling a situation in which there are several large groups of electricity producers and in each group there is a major player who can influence the cost function of the producers. Each producer is considered as a minor player in the sense that she has a small impact on the global behavior of the system. The populations interact with a network structure, meaning that each group interact with some but not necessarily all other groups. We can think of such a model as representing electricity producers in different countries, and in each country the government can influence the producers’ cost function by setting taxes. The government itself tries to optimize a criterion. In the present work, we consider that the government is in a Nash equilibrium with the producers. We want to stress that this is different than Stackelberg equilibrium between major players and mean field populations in which the major players optimize over their cost functions after anticipating the effect of their controls on the behavior of minor players.

Our main contributions are three fold. First we introduce a stylized model for carbon emissions with multiple mean field populations and one major player per population, with network interactions. Second, relying on Pontryagin’s maximum principle, we provide analytical solutions for the minor players' controls and the major players' controls at equilibrium. Last, we illustrate numerically the sensitivity of the model with respect to key parameters.

The rest of the paper is organized as follows. In Section~\ref{sec:model} we introduce the general model, starting from the finite player setting before presenting the mean field game. In Section~\ref{sec:carbon}, we present a model for carbon emission regulation as an example of the general model, and we solve it explicitly. In Section~\ref{sec:algo-num}, we present numerical results that illustrate the sensitivity of the model to its parameters. Last, in Section~\ref{sec:conclusion-future}, we conclude and discuss future works.

\section{Model}
\label{sec:model}

\subsection{Model With Finitely Many Players}
\label{subsec:model_finite}
First, we start by introducing the finite player model which motivates the mean field model. Let $T>0$ be the finite time horizon, let $M\geq 2$ be the number of populations, and $N_i>0$ be the number of minor players in each population $i$. There is a major player in each population $i$, for all $i\in \{1,\dots, M\}=\llbracket M\rrbracket$.  Furthermore, let $\bsG$ be the connection matrix that gives the connection strength between populations, i.e., $\bsG_{k,\ell}\in[0,1]$ is the connection level between population $k$ and $\ell$. 

Minor player $j$ in population $i$ controls her state $(X^{i,j}_t)_{t\in[0,T]}$ by choosing her control $\alpha^{i,j} = (\alpha_t^{i,j})_{t\in[0,T]}$. We denote by $\bsalpha =(\alpha^{i,j})_{i=1,\dots,M,j=1,\dots,N_i}$ the control profile for all the minor agents. Furthermore, we will denote by $\gamma^i$ the control of the major player $i$. We denote by $\bsgamma = (\gamma^i)_{i=1,\dots,M}$ the control profile for the major players. Given $(\bsalpha,\bsgamma)$, the state of the minor player $j$ in population $i$ follows the dynamics: $X^{i,j}_0 = x^{i,j}_0 \sim \mu_0^i$ and
\begin{equation*}
    dX^{i,j}_t = b_{\bsG}(t, X_t^{i,j}, \alpha_t^{i,j}, \bar{X}^N_t, \bar{\alpha}^N_t, \gamma^i)dt + \sigma^i dW_t^{i,j},
\end{equation*} 
where $\mu_0^i$ is the initial state distribution of players in population $i$ and the initial positions are independent, 
$b:[0,T] \times \RR \times \RR \times \RR^M \times \RR^M \times \RR \to \RR$ is the drift, $\sigma^i\in \RR$ is the volatility and $W^{i,j}$ is a Brownian motion representing some idiosyncratic noise, independent of the Brownian motions for the other agents. Here, $\bar{X}^N_t = (\frac{1}{N_1}\sum_{j=1}^{N_1} X^{1,j}_t,\dots \frac{1}{N_M}\sum_{j=1}^{N_M} X^{M,j}_t)$ is the average states, and $\bar{\alpha}^N_t = (\frac{1}{N_1}\sum_{j=1}^{N_1} \alpha^{1,j}_t,\dots \frac{1}{N_M}\sum_{j=1}^{N_M} \alpha^{M,j}_t)$ is the average actions in the populations.

Then, the cost incurred to player $i$ is:
\begin{equation*}
    \begin{aligned}
    &J^N(\bsalpha^{i,j}; \bsalpha^{-(i,j)}, \gamma^i)=\\ 
    &\EE\Big[\int_0^T f_{\bsG}(t, X_t^{i,j}, \alpha_t^{i,j}, \bar{X}^N_t, \bar{\alpha}^N_t, \gamma^i)dt
    + g_{\bsG}(X_T^{i,j}, \bar{X}_T, \gamma^i)\Big]
    \end{aligned}
\end{equation*}
where $\bsalpha^{-(i,j)}$ denotes the control of every minor agent except agent $(i,j)$. Here, $f_{\bsG} : [0,T] \times \RR \times \RR \times \RR^M \times \RR^M \times \RR \to \RR$ and $g_{\bsG} : \RR \times \RR^M \times \RR \to \RR$ are the running cost and terminal cost functions for minor players, respectively. They depend on the network structure $\bsG$ between populations.

We now introduce the notion of Nash equilibrium when only the minor players are optimizing. However, in the sequel, we will consider a model in which the major players can also optimize their cost.
\begin{definition}
\label{def:minor-NE-finite}
Given a control profile $\bsgamma = (\bsgamma^{i})_{i}$ for the major players, a minor players' control profile $\hat\bsalpha = (\hat\bsalpha^{i,j})_{i,j}$ 
is called a multi-population minor Nash equilibrium if for every $i$ and every $j$, $\hat\alpha^{i,j}$ minimizes the minor player cost: $\alpha^{i,j} \mapsto J^N(\alpha^{i,j}; \hat\bsalpha^{-(i,j)}, \gamma^i)$. 
\end{definition}

The major player in population $i$ pays the following cost when she chooses a control level $\gamma^i$:
\begin{equation*}
\begin{aligned}
    &J^{N,0}(\gamma^i; \bsgamma^{-i}, \bsalpha)=\\
    &\EE\left[\int_0^T f^{0}_{\bsG}(t, \gamma^i, \gamma^{-i}, \bar{X}^N_t, \bar{\alpha}_t)dt + g^{0}_{\bsG}(\gamma^i, \gamma^{-i}, \bar{X}^N_t)\right]
\end{aligned}
\end{equation*}
where $\bsgamma^{-i}$ is the vector of control of other major players i.e., $\bsgamma^{-i} = (\gamma^1, \dots, \gamma^{i-1}, \gamma^{i+1}, \dots, \gamma^M)$. Here,  $f^{0}_{\bsG}: [0,T]\times \RR \times \RR^{M-1} \times \RR^M \times \RR^M \to \RR$ and $g^{0}_{\bsG} : \RR \times \RR^{M-1} \times \RR^M$ are the running cost and terminal cost functions of major players, respectively. In this model, it would be possible to let $f$, $g$, $f^{0}$ and $g^{0}$ depend on $i$, but we omit this to alleviate the presentation. Furthermore, in this model the major players only need to choose one action. The extension to major players with states and time-dependent controls is left for future work. From our motivated application point of view (i.e., choosing the carbon tax levels) in this study, having a constant control for major players is a natural choice since, these levels will be fixed over the fixed time interval.

We can see that in the model minor player $j$ in population $i$ interacts with the minor agents in her population as well as the minor agents in other populations through the population mean of the states and controls ($\bar{\bsX}$ and $\bar{\bsalpha}$). She also interacts with the major player $i$ through the major player's control ($\gamma^i$). On the other hand, major player $i$ interacts with the other major players through their controls ($\bsgamma^{-i}$) and with the minor players through their population mean of states and controls.

\begin{definition}
\label{def:major-minor-mfne} 
A minor and major player control tuple $(\hat\bsgamma, \hat\bsalpha)=(\hat\bsgamma^{i}, \hat\bsalpha^{i,j})_{i,j}$ 
is called a multi-population major minor Nash equilibrium if
\begin{itemize}
    \item For every $i$ and every $j$, $\hat\alpha^{i,j}$ minimizes the minor player cost: $\alpha^{i,j} \mapsto J^N(\alpha^{i,j}; \hat\bsalpha^{-(i,j)}, \hat\gamma^i)$,
    \item For every $i$, $\hat\gamma^i$ minimizes the major player cost:  $\gamma^i \mapsto J^{N,0}(\gamma^i; \hat\bsgamma^{-i}, \hat\bsalpha)$.
\end{itemize}
\end{definition}

\subsection{Mean Field Game Model}
\label{subsec:model_meanfield}
When the number of minor players is very large, each minor player has a very small impact on the other players since the interactions are only through the average state and action. Keeping track of all pairwise interactions and computing an exact Nash equilibrium would be computationally prohibitive, so we rely on a mean-field approximation in order to have a tractable formulation. We assume that the number of people in each of the populations goes to infinity, and focus on the limiting mean field models for each of the populations. In the limiting regime, the behavior of population $i$ is summarized by the behavior of one representative minor player in this population. Next, we consider that $N\rightarrow \infty$.

A representative minor player in population $i$ controls her state $(X^{i}_t)_{t\in[0,T]}$ by choosing her control $\alpha^{i} = (\alpha_t^{i})_{t\in[0,T]}$. We denote by $\bsalpha = (\alpha^{i})_{i=1,\dots,M}$ the control profile for the minor agents. We continue to denote by $\gamma^i$ the control of the major player $i$ and by $\bsgamma = (\gamma^i)_{i}$ the control profile for all the major players. Given $(\bsalpha,\bsgamma)$, the state of a minor player in population $i$ follows the dynamics:
\begin{equation}
    dX^{i}_t = b_{\bsG}(t, X_t^{i}, \alpha_t^{i}, \bar{X}_t, \bar{\alpha}_t, \gamma^i)dt + \sigma^i dW_t^{i},
\end{equation}
where $b$ and $\sigma^i$ are as before, and $W_t^{i}$ is a Brownian motion representing some idiosyncratic noise, independent of the Brownian motions for the other agents. Here, $\bar{X}_t = (\EE[X^1_t], \dots, \EE[X^M_t])$ is the average states, and $\bar{\alpha}_t = (\EE[\alpha^1_t], \dots, \EE[\alpha^M_t])$ is the average actions in the populations. Then, the cost incurred to a minor player in population $i$ is: 
\begin{equation}
    \begin{aligned}
    &J(\alpha^{i}; \bar{\boldsymbol{X}}, \bar{\bsalpha}, \gamma^i)=\\ 
    &\EE\Big[\int_0^T f_{\bsG}(t, X_t^{i}, \alpha_t^{i}, \bar{X}_t, \bar{\alpha}_t, \gamma^i)dt
    + g_{\bsG}(X_T^{i}, \bar{X}_T, \gamma^i)\Big]
    \end{aligned}
\end{equation}
where $\bsalpha^{-i} = (\alpha^{1}_t,  \dots,\alpha^{i-1}_t,\alpha^{i+1}_t, \dots,\alpha^{N_M}_t)_{t\in[0,T]}$.  Here, $f_{\bsG}$ and $g_{\bsG}$ are defined as before and depend on the network structure between populations.
We stress that, in this model, the interactions involve both the state distribution and the control distribution (through their means). This type of setting is sometimes referred to as extended MFG~\cite{gomes2016extended} or MFG of controls~\cite{cardaliaguet2018mean,kobeissi2022classical}. We can now introduce the counterpart in the mean field regime of Definition~\ref{def:minor-NE-finite}.
\begin{definition} 
\label{def:minor-NE-meanfield}
Given a control profile $\bsgamma = (\gamma^{i})_{i}$ for the major players, a minor players' control profile $\hat\bsalpha = (\hat\alpha^{i})_{i}$ 
is called a \textit{multi-population mean field Nash equilibrium} if for every $i$, $\hat\alpha^{i}$ minimizes the minor player cost: $\alpha^{i} \mapsto J(\alpha^{i}; \bar{\boldsymbol{X}}, \bar{\bsalpha}, \gamma^i)$. 
\end{definition}

The major player in population $i$ pays the following cost when she chooses a control level $\gamma^i$:
\begin{equation}
\begin{aligned}
    &J^{0}(\gamma^i; \bsgamma^{-i}, \bar{\boldsymbol{X}}, \bar{\bsalpha})=\\
    &\int_0^T f^{0}_{\bsG}(t, \gamma^i, \gamma^{-i}, \bar{X}_t, \bar{\alpha}_t)dt + g^{0}_{\bsG}(\gamma^i, \gamma^{-i}, \bar{X}_T)
\end{aligned}
\end{equation}
where $\bsgamma^{-i}$ is the vector of controls of other major players i.e., $\bsgamma^{-i} = (\gamma^1, \dots, \gamma^{i-1}, \gamma^{i+1}, \dots, \gamma^M)$. Here,  $f^{0}_{\bsG}$ and $g^{0}_{\bsG}$ are defined as before.

\begin{definition} A minor and major player control tuple $(\hat\bsgamma, \hat\bsalpha)=(\hat\bsgamma^{i}, \hat\bsalpha^{i})_{i}$ is called a \textit{multi-population major minor mean field Nash equilibrium} (M4FNE) if
\begin{itemize}
    \item For every $i$, $\hat\alpha^{i}$ minimizes the minor player cost: $\alpha^{i} \mapsto J(\alpha^{i}; \bar{\boldsymbol{X}}, \bar{\bsalpha}, \gamma^i)$
    \item For every $i$, $\hat\gamma^i$ minimizes the major player cost:  $\gamma^i \mapsto J^{0}(\gamma^i; \bsgamma^{-i}, \bar{\boldsymbol{X}}, \bar{\bsalpha})$.
\end{itemize}
\end{definition}

\section{Application to the control of carbon emission levels}
\label{sec:carbon}

In this section, we introduce our motivating example for finding the equilibrium carbon tax levels for interacting countries, discuss its solution characterization and give the theoretical results. All the proofs are presented in the appendix.

\subsection{Mean Field Game Model}
In our model, we have electricity producers as the minor players and the tax regulators as the major players in populations $1, \dots, M$. The state of the representative electricity producer in population $i$ at time $t$ is the cumulative production up until time $t$ and their control is to how much nonrenewable energy they use at time $t$. 

The dynamics of the state is given as:
\begin{equation}
    dX_t^i = \eta^i \alpha_t^i dt + \sigma^i dW^i_t,\quad X^i_0 = x^i_0 \sim \mu_0^i,
\end{equation}
where $\eta^i>0$ is a constant that gives the efficiency of one unit of nonrenewable energy in the production of electricity.  Intuitively, the dynamics states that the cumulative production until time $t$ is equal to the multiplication of efficiency of the nonrenewable energy and the total nonrenewable energy used until time $t$. In order to account for the randomness in the production levels that comes for example from facility conditions, we add an idiosyncratic Brownian motion.

The representative electricity producer in population $i$ incurs the following cost when she chooses to use the control $(\alpha_t^i)_{t\in[0,T]}$:
\begin{equation}
    \begin{aligned}
        J(\alpha^{i}; \bar{\boldsymbol{X}}, \bar{\bsalpha}, \gamma^i) =& \EE\Big[\int_0^T \Big( \gamma^i ({\alpha}^i_t)^2
        + 
        \sum_{i'=1}^{M} \textbf{G}_{i,i'} (\bar \alpha_t^{i'})^2 
        \Big) dt
        \\
        &\hskip2mm
        + 
        \kappa^i\sum_{i'=1}^{M} \textbf{G}_{i,i'}({X}^{i}_T-\bar{X}_T^{i'})^2 -\delta^i X^{i}_T\Big],
    \end{aligned}
\end{equation} 
where $\kappa^i$, $\delta^i>0$ are constants that balance out the importance of the cost terms. Here, the first term in the integral is the carbon tax cost that is paid depending on the nonrenewable energy usage. We assume that the producers are getting taxed more at higher levels of nonrenewable energy usage; therefore, the control of the minor is squared. The second term in the integral can be thought as the cost of pollution where the effect of the pollution from other countries is also included. The effects from different countries are weighted according to the connection matrix matrix $\bsG$. The third term compares the production level of the representative minor player in population $i$ to the average production levels in other populations. If the production level of the representative minor player in population $i$ is lower than the average production levels of their neighbours, it behaves as an opportunity cost of not producing more. If the production level of representative minor player in population $i$ is higher than the average production levels of their neighbours, it behaves as a reputation cost for polluting more. Finally, the last term accounts from net earnings from the cumulative production at the terminal time $T$.

The major player $i$ controls the tax level and incurs the following cost when she chooses to use the control $\gamma^i$:
\begin{equation}
  \begin{aligned}  
    J^{0}(\gamma^i; \bsgamma^{-i}, \bar{\boldsymbol{X}}, \bar{\bsalpha})
    = &\kappa_c J(\alpha^{i}; \bar{\boldsymbol{X}}, \bar{\bsalpha}, \gamma^i)  
    - \int_0^T \gamma^i\EE[(\alpha_t^i)^2] dt
    \\
    &\hskip2mm +\sum_{i'=1}^{M} \textbf{G}_{i,i'} (\gamma^i-\gamma^{i'})^2
    + \kappa_g (\gamma^i)^2,
\end{aligned}
\end{equation}
where $\kappa_c$ and $\kappa_g>0$ are constants that balance out the weights of different cost terms. The major player gets negatively affected when the minor players in their population is incurring high costs. The first term accounts for this fact. The second term is the average carbon tax revenue. The third term compares the carbon tax level in population $i$ to the carbon tax levels in other populations. If the carbon tax level in population $i$ is lower than their neighbours, it behaves as an opportunity cost of not taxing more and if the carbon tax level in population $i$ is higher than their neighbours, it behaves as a reputation cost. The last term is the reputation cost of charging high levels of tax.

\subsection{Theoretical Analysis}
We now present a theoretical analysis of the model. We first note that, for a given profile $\bsgamma = (\gamma^i)_{i}$ for the major players, the game for the minor players is a multi-population MFG for which the solution reduces to a system of forward-backward stochastic differential equations in which the backward variable represents the derivative for the value function for one representative player. 
\begin{theorem}
\label{the:minors-FBSDE}
A control profile $\hat{\bsalpha}$ is an equilibrium control profile  (see Definition~\ref{def:minor-NE-meanfield}) if and only if:
    \begin{equation}
    \label{eq:BR-formula}
        \hat \alpha^{i}_t
        = 
        - \frac{\eta^i}{2 \gamma^i} Y^i_t, \quad \forall i \in\llbracket M \rrbracket
    \end{equation}
    where $(\bsX, \bsY, \bsZ)=(X_t^i, Y_t^i, Z_t^i)_{i\in \llbracket M \rrbracket,\ t\in [0,T]}$ solve the following forward-backward stochastic differential equation (FBSDE) system:
    \begin{equation}
    \label{eq:FBSDE_minor}
    \begin{cases}
        dX^i_t = -\dfrac{{|\eta^i|^2}}{2\gamma^i}Y^i_tdt +\sigma^i dW^i_t, 
        \hfill X^i_0 = x^i_0,\ \forall i \in \llbracket M \rrbracket, &
        \\[2mm]
        dY^i_t = Z^i_t d W^i_t, \ &
        \\[1mm]
        Y^i_T = -\delta^i + 2\kappa^i (\bsG[i,:] \mathds{1}_M  X^i_T-\bsG[i,:]\bar{X}_T),\hfill  \forall i \in \llbracket M \rrbracket,&
    \end{cases}
    \end{equation}
    where $\mathds{1}_M$ is the (column) vector of ones with length $M$, $\bsG[i,:]$ is the $i^{\rm th}$ row of connection matrix $\bsG$. 
\end{theorem}
The proof is reminiscent of the proof of the stochastic Pontryagin's maximum principle, see~\cite[Section 4.5]{carmona2018probabilistic} in a fully homogeneous MFG and~\cite[Section 7.1.1]{carmona2018probabilistic} for multi-population MFG. We adapt the proof technique to our setting, which involves interactions through the controls' distribution instead of the state distribution.

Using the linear quadratic structure, we can reduce the computation of its solution to the solution of a forward-backward ODE system, as we explain next. 
\begin{proposition}
\label{prop:minors-ode}
Given $\bsgamma = (\gamma^1, \dots, \gamma^M)$, assume there exists an $\RR^{M} \times \RR^{M \times M} \times \RR^{M}$-valued function $t \mapsto (A_t, B_t, \bar{X}_t)$ with $\bar{X}_t = (\bar{X}_t^1, \dots, \bar{X}_t^M)$ solving the following forward-backward ordinary differential equation (FBODE) system:
    \begin{align}
        &{\dot A^i_t}  =  B^i_t \bsG[i,:]\mathds{1}_M \dfrac{|\eta^i|^2A^i_t}{2{\gamma^i}}\nonumber\\
        &-B^i_t \sum_{j=1}^M \dfrac{|\eta^j|^2\bsG[i,j]}{2{\gamma^j}}\Big[A^j_t+B^j_t\left(\bsG[j,:](\mathds{1}_M \bar{X}^j_t - \bar{X}_t)\right) \Big], \label{eq:FBODE_A} \\
        &\dot B^i_t  = \dfrac{{|\eta^i|^2}(B^i_t)^2}{2{\gamma^i}}{\bsG}[i,:]\mathds{1}_M, \label{eq:FBODE_B}
       \\
       &\dot{\bar{X}}^i_t =  -\dfrac{{|\eta^i|^2}}{2\gamma^i}\Big( A_t^i + B^i_t \big[\bsG[i,:]\mathds{1}_M\bar{X}^i_t - \bsG[i,:]\bar{X}_t\big] \Big), \label{eq:FBODE_BarX}\\
       &A^i_T =-\delta^i, \quad B_T^i = 2\kappa^i, \quad \bar{X}_0^i = \bar x_0^i \label{eq:FBODE_boundarycond}
    \end{align}
for all $i\in\llbracket M\rrbracket$. Then 
\begin{equation}
    \label{eq:equilibrium-ctrl-minors-only}
\hat \alpha_t^i =  - \frac{{\eta^i}}{2 \gamma^i} \big(A_t^i + B_t^i(\bsG[i,:]\mathds{1}_MX_t^i - \bsG[i,:]\bar{X}_t)\big)
\end{equation}
is the equilibrium control of minor player $i$ for all $i\in \llbracket M \rrbracket$. 
\end{proposition}
The proof consists of introducing an ansatz for $Y^i_t$ as an affine function of $X^i_t$ in the form $Y^i_t = A^i_t  + B^i_t (\bsG[i,:] \mathds{1}_M X^i_t - \bsG[i,:]\bar{X}_t)$, which is consistent with the fact that $Y^i_t$ represents the derivative of the value function of a representative player in population $i$ and, in linear-quadratic models, we expect the value function to be quadratic.

\begin{remark}
The average control of minor players in population $i$ is given by
\begin{equation}
\label{eq:minor_control_avg}
    \bar{\hat{\alpha}}^i_t = - \frac{{\eta^i}}{2 \gamma^i} \big(A_t^i + B_t^i(\bsG[i,:]\mathds{1}_M\bar{X}_t^i - \bsG[i,:]\bar{X}_t)\big).
\end{equation}
\end{remark}

Next, we prove existence and uniqueness of an equilibrium for the game between minor players. 

\begin{theorem}
\label{thm:existence-uniqueness}
Given any major players' control profile $\bsgamma = (\gamma^i)_i$ with $\gamma^i>0$ for all $i=1,\dots,M$, there exists $T_0$ such that for all $T<T_0$, there exists a unique mean field Nash equilibrium between the minor players.  
\end{theorem}

The proof relies on an application of Banach fixed point theorem after proving upper bounds on the norm of $\bar{X}^i$, $i=1,\dots,M$.

We then derive explicit expressions for the major players optimal controls in terms of the other players' controls. 
\begin{theorem}
\label{thm:majors-controls-formula}
    Given the controls of representative minor players, $\bsalpha$, and  other major players $\bsgamma^{-i}$, the best response of the major player $i$ is given by: 
    \begin{equation}
    \label{eq:major-best-response}
        {\gamma}^{i*} = \dfrac{\tau^i - 2\sum_{j=1, {j\neq i}}^{M}\bsG_{i,j}\gamma^j}{{-}2\kappa_g-2\bsG[i,:]\mathds{1}_M + 2\bsG[i,i]},
    \end{equation}
    where $\tau^i = (\kappa_c-1) \int_0^T \mathbb{E}[(\alpha_t^i)^2]dt$. Given the controls of representative minor players, $\bsalpha$, the equilibrium control profile for the major players is:
    \begin{equation}
    \label{eq:major-eq-given-minor}
    \hat{\bsgamma} = - \frac{1}{2} \big[\diag(\bsG \mathds{1}_M + \kappa_g \mathds{1}_M) -\bsG \big]^{-1}\boldsymbol{\tau},
    \end{equation}
    provided the matrix $\big[\diag(\bsG \mathds{1}_M + \kappa_g \mathds{1}_M) -\bsG \big]$ is invertible. 
\end{theorem}

\section{Algorithm and Numerical Results} 
\label{sec:algo-num}

\subsection{Algorithm}

In order to compute the global Nash equilibrium between all the players (M4FNE), we can rely on Proposition~\ref{prop:minors-ode} and Theorem~\ref{thm:majors-controls-formula} respectively for the minor players and the major players. Note that the two expressions~\eqref{eq:equilibrium-ctrl-minors-only} and~\eqref{eq:major-eq-given-minor} are coupled, so we cannot solve one before the other. Furthermore, the former requires solving a coupled forward-backward ODE system. To compute the equilibrium, we use a fixed point algorithm: we iteratively update the minors' controls and the majors' controls. To solve numerically the FBODE system in Proposition~\ref{prop:minors-ode}, we discretize time using an Euler scheme. Here the equation for $B^i$ is uncoupled and can be explicitly solved. However, since the forward equation for $A^i$ and the backward equation for $\bar{X}^i$ are coupled, here again we cannot solve one before the other. For this step, we use a subroutine which iteratively updates each variable until convergence. The pseudo-code of the algorithm is summarized in Algorithm~\ref{algo:SGD-SMFG}. We want to stress that this algorithm has two nested fixed point iterations where one of them is nested in the other. If one is interested in finding the multi-population mean field Nash equilibrium when the major players' controls are given exogenously, the FBODE system in Proposition~\ref{prop:minors-ode} needs to be solved for given controls of the major players. In this setup, only the inner fixed point iteration that solves coupled FBODE should be used.

{
\begin{algorithm}
\caption{\small Multi-population Major Minor Mean Field Nash Equilibrium (M4FNE) Algorithm\label{algo:SGD-SMFG}}

\small
\textbf{Input:} Parameters of the minor players' model: $(\delta^i)_{i \in \llbracket M \rrbracket}, \kappa, \eta$; connection matrix: $\bsG$; parameters of the major players' model: $\kappa_c, \kappa_g$; tolerance: $\epsilon$; time horizon: $T$; time increments: $\Delta t$.\footnote{We will denote the time steps $\{0, 1, \dots, n_T\}$ by $\llbracket n_T \rrbracket$ where $T = n_T \Delta t$.}

\textbf{Output:} Equilibrium carbon tax levels, $\bsgamma$; Equilibrium mean nonrenewable energy usage of producers, $\bar{\bsalpha}$; Equilibrium average cumulative production, $\bar{\bsX}$.
\vskip1mm

\begin{algorithmic}[1]

\STATE{Initialize $\bsgamma^{(0)} = (\gamma^{(0),1}, \dots, \gamma^{(0),M})$ and $\bar{\bsalpha}^{(0)} =(\bar{\alpha}_t^{(0),i})_{i \in \llbracket M \rrbracket, t=n_t\Delta t, n_t \in \llbracket n_T\rrbracket}$.} \vskip2mm

\WHILE{$||\bsgamma^{(j)}-\bsgamma^{(j-1)}||>\epsilon$ or $||\bar{\bsalpha}^{(j)}-\bar{\bsalpha}^{(j-1)}||>\epsilon$}
    \STATE{Solve the backward ODE in \eqref{eq:FBODE_B} given $\bsgamma^{(j)}$ for $\bsB^{(j+1)}$ $=(B_t^{(j+1),i})_{i \in \llbracket M \rrbracket, t\in\llbracket n_T\rrbracket\Delta t}$ }  \vskip1mm
    \STATE{Solve the coupled forward-backward ODEs in \eqref{eq:FBODE_A} and \\\eqref{eq:FBODE_BarX} given $\bsgamma^{(j)}$ and $\bsB^{(j+1)}$ for $\bsA^{(j+1)}=(A_t^{(j+1),i})_{i \in \llbracket M \rrbracket, t\in\llbracket n_T\rrbracket\Delta t}$ and $\bar{\bsX}^{(j+1)}=(\bar{X}_t^{(j+1),i})_{i \in \llbracket M \rrbracket, t\in\llbracket n_T\rrbracket\Delta t}$}\vskip1mm
    \STATE{Compute $\bar{\bsalpha}^{(j+1)}$ given $\bsgamma^{(j)}, \bsA^{(j+1)}, \bsB^{(j+1)},$ and\\ $\bar{\bsX}^{(j+1)}$ using \eqref{eq:minor_control_avg}.}\vskip1mm
    \STATE{Compute best response of the major players $\bsgamma^{(j+1)}$ \\given $\bar{\bsalpha}^{(j+1)}$.}
\ENDWHILE

\RETURN $\bsgamma^{(j+1)}, \bar{\bsalpha}^{(j+1)},\bar{\bsX}^{(j+1)} $
\end{algorithmic}
\end{algorithm}}

\subsection{Numerical Results}

To illustrate our results, we provide some numerical examples. In our numerical example, we look at a setup where we have 3 populations and 3 major players. For the following experiments, we take the initial means to be $\bar{x}_0^1=25, \bar{x}_0^2=20$ and $\bar{x}_0^3=15$ respectively for populations $1$ and $2$, and $\eta^i=1, \kappa^i=0.1$ for all populations. 

\begin{figure}[t]
     \centering
     \begin{subfigure}[b]{0.45\textwidth}
         \centering
         \includegraphics[width=\textwidth]{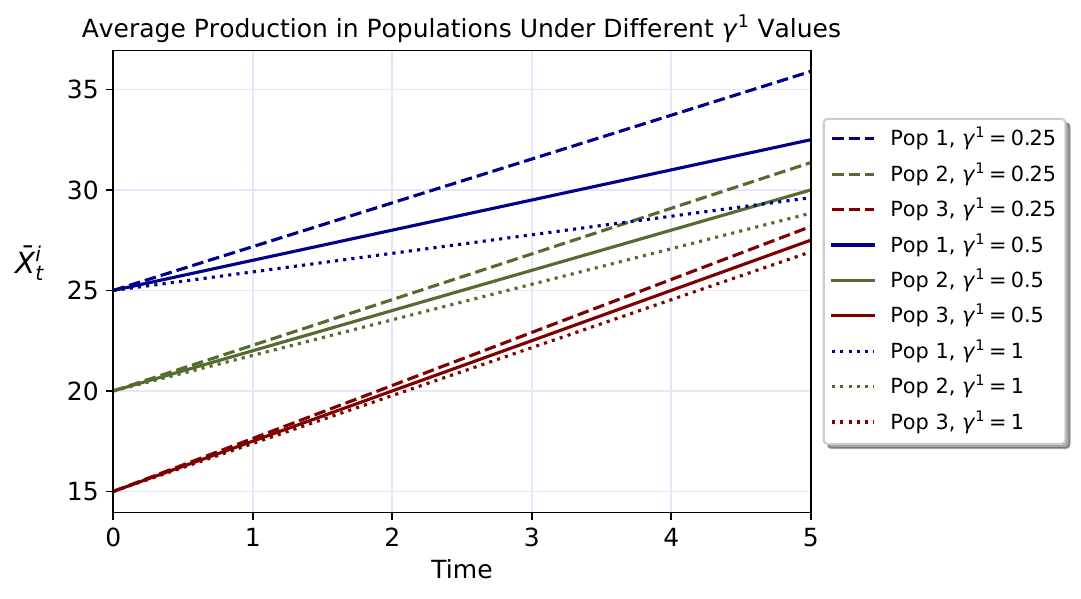}
     \end{subfigure}
     \begin{subfigure}[b]{0.45\textwidth}
         \centering
         \includegraphics[width=\textwidth]{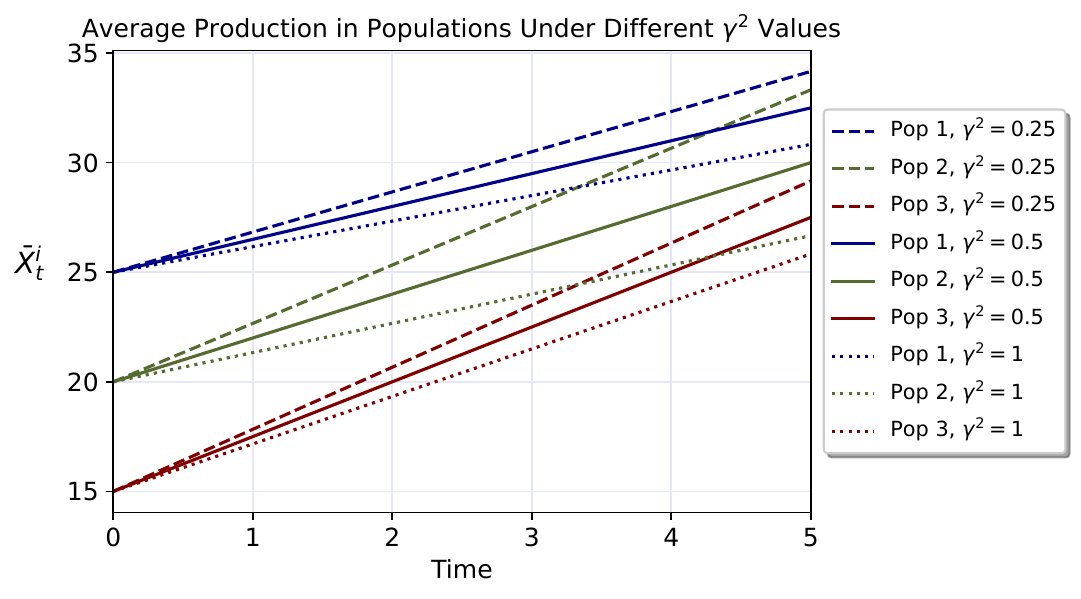}
     \end{subfigure}
        \caption{\small The effect of different carbon tax levels in population 1 (top) and in population 2 (bottom) on the average production levels in the populations.}
        \label{fig:minor_gamma_effect}
\end{figure}

\begin{figure}[t]
     \centering
     \begin{subfigure}[b]{0.45\textwidth}
         \centering
         \includegraphics[width=\textwidth]{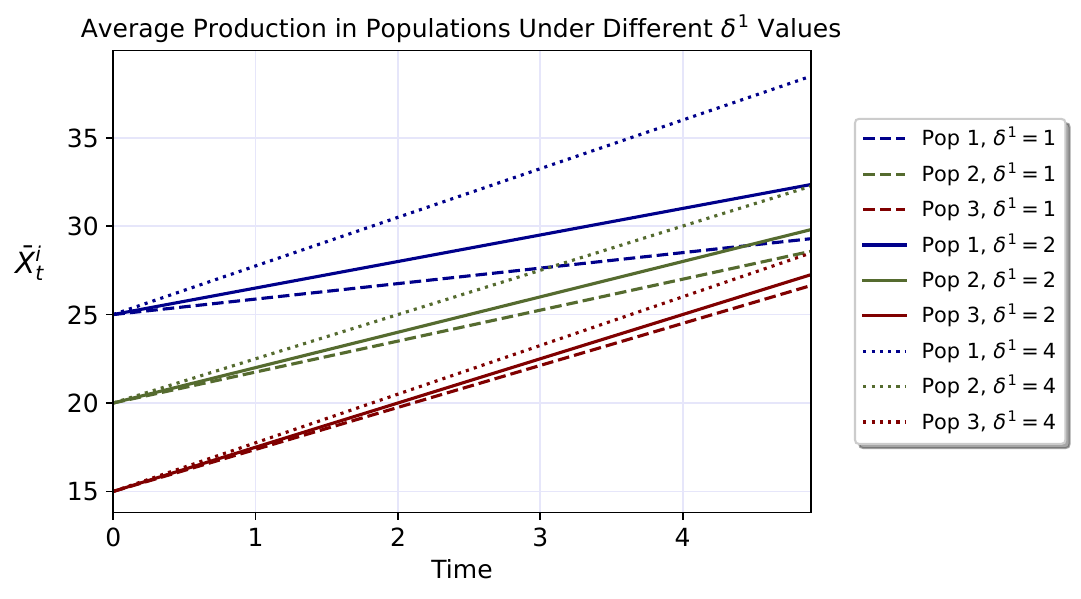}
     \end{subfigure}
     \begin{subfigure}[b]{0.45\textwidth}
         \centering
         \includegraphics[width=\textwidth]{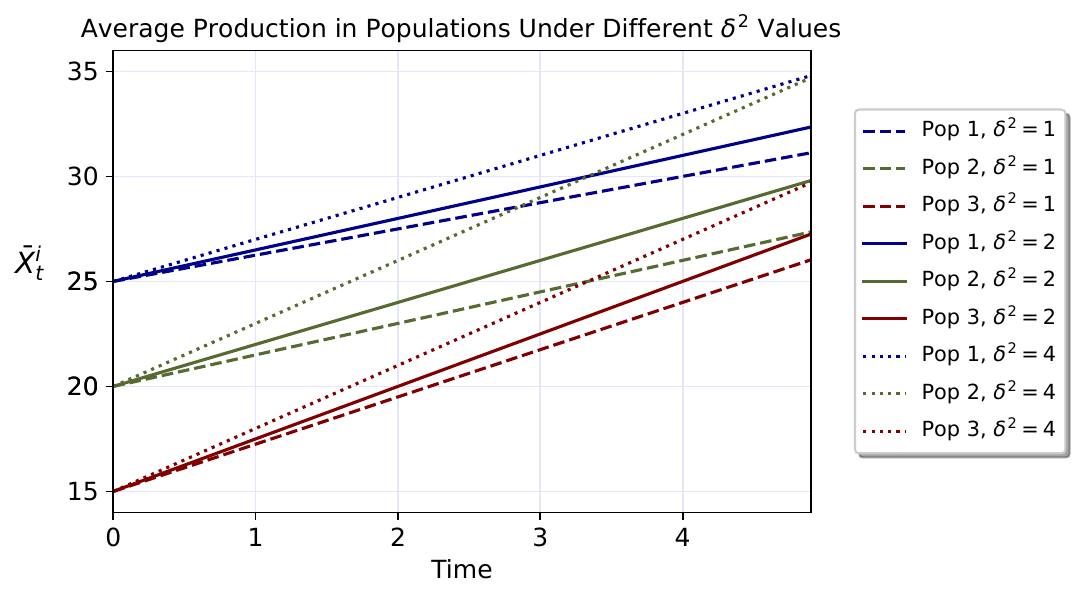}
     \end{subfigure}
        \caption{\small The effect of different coefficient of net earnings in population 1 (top) and in population 2 (bottom) on the average production levels in the populations.}
        \label{fig:minor_delta_effect}
\end{figure}

\begin{figure}[b]
     \centering
     \begin{subfigure}[b]{0.45\textwidth}
         \centering
         \includegraphics[width=\textwidth]{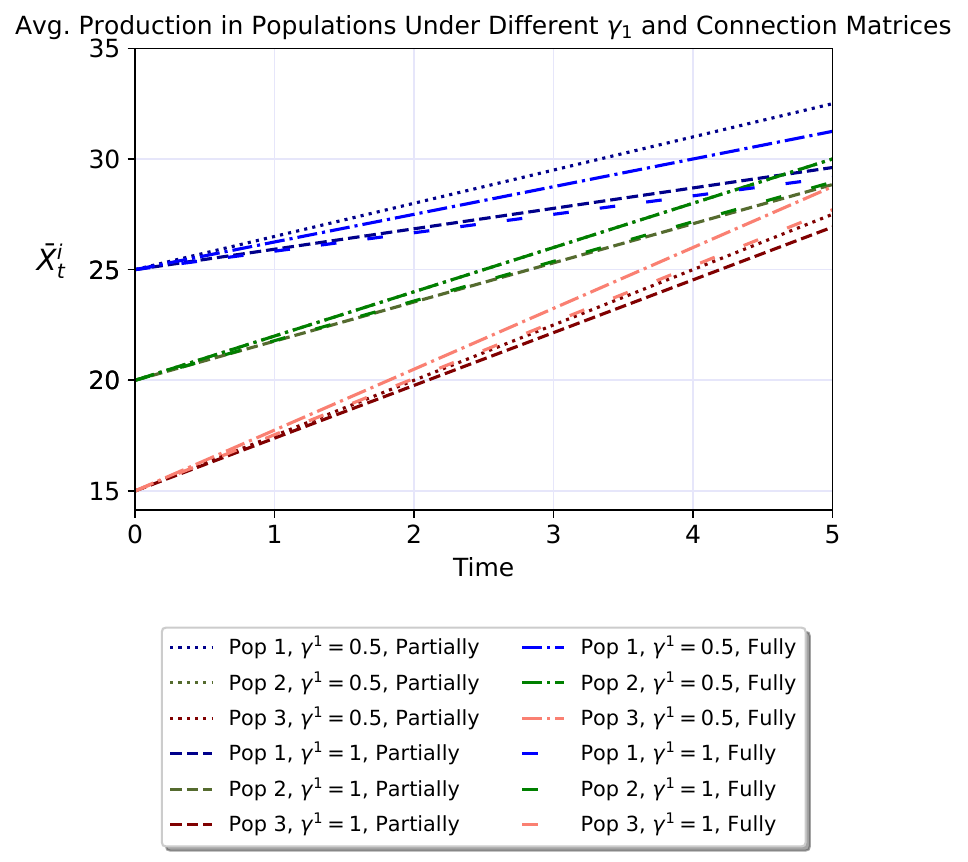}
     \end{subfigure}
        \caption{\small The effect of different carbon tax levels in population 1 on the average population levels when there are different connection setups between the populations.}
        \label{fig:minor_connection_gamma_effect}
\end{figure}

We first analyze the effect of the coefficient of net earnings ($\delta^i$) and carbon tax levels ($\gamma^i$) on the multi-population mean field Nash equilibrium when the carbon tax levels are set exogeneously. For these experiments, we are assuming that population 1 and population 2, and population 2 and population 3 are fully connected to each other, while the population 1 and population 3 are not connected to each other (i.e., we take $\bsG_{i,j}=0$ for $(i,j) \in\{(1,3),(3,1)\} $ and $\bsG_{i,j}=1$, otherwise). In Fig.~\ref{fig:minor_gamma_effect}, we can see the effect of different carbon tax levels in population 1 (top figure) and in population 2. As expected, when the carbon tax in one population is increased, the production in that population decreases. This further affects the populations other populations interacting. In Fig.~\ref{fig:minor_delta_effect}, we can see that increasing the coefficients of net earnings from production in population 1 (top figure) and in population 2 (bottom figure) results in higher production levels in the populations the coefficient is increased. We can also see that the production is increased (but less significantly) in other countries because of the comparison cost among the populations. In both experiments, we can see a change in population 2 affects population 3 further than a change in population 1 because population 1 and 3 are only indirectly connected through the population 2. Furthermore, we check the effect of connections, where we compare the current \textit{partially connected} populations setup to the \textit{fully connected} populations setup (i.e., we take $\bsG_{i,j}=1$ for all $i, j$ in the fully connected setup). In Fig.~\ref{fig:minor_connection_gamma_effect}, we can see that when populations are fully connected an increase the carbon tax levels in population 1 affects the production in population 3 more than the partially connected setup.

\begin{figure}
     \centering
     \begin{subfigure}[b]{0.45\textwidth}
         \centering
         \includegraphics[width=\textwidth]{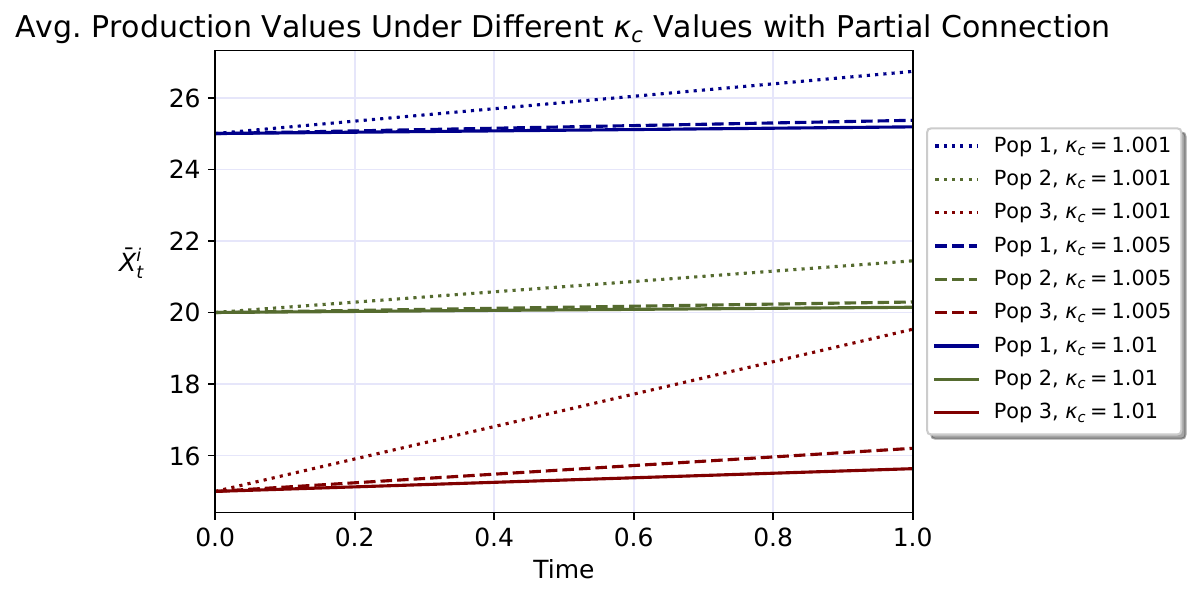}
     \end{subfigure}
        \begin{subfigure}[b]{0.45\textwidth}
         \centering
         \includegraphics[width=\textwidth]{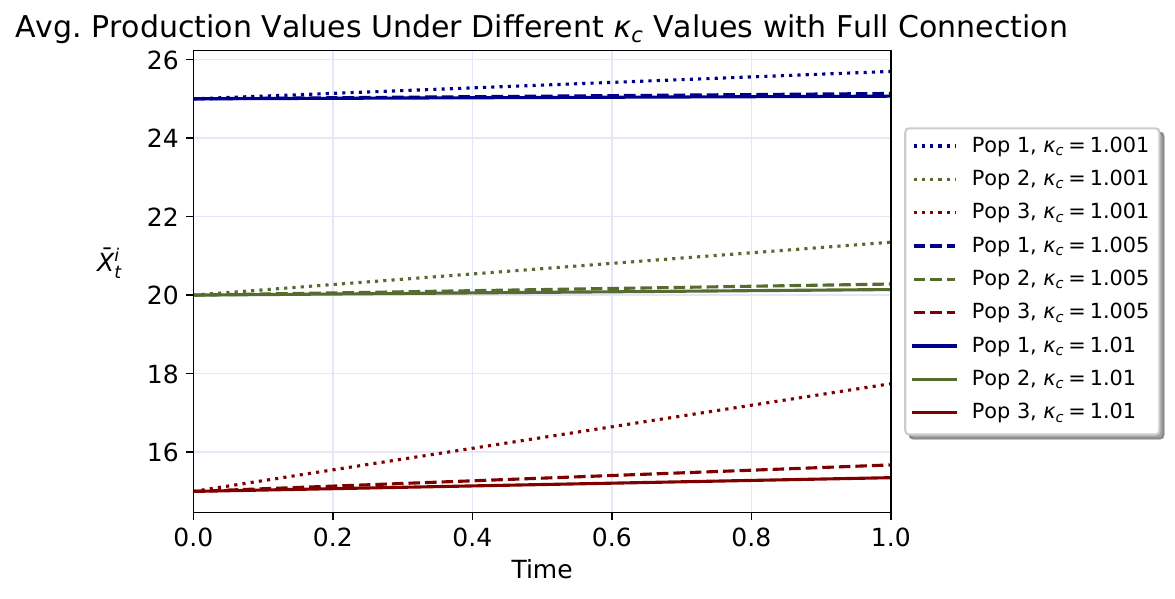}
     \end{subfigure}
        \caption{\small The effect of different coefficient of importance major players' give to the cost of minor players in their population on the average production levels in the populations when populations are partially connected (top) or fully connected (bottom).}
        \label{fig:major_connection_effect}
\end{figure}

Later, we focus on experiments where we find the M4FNE and check the effect of the coefficient of importance major players' give to the cost of the minor players ($\kappa_c$). The M4FNE carbon tax levels can be seen in Table~\ref{tab:MajMin-NE_tax} under different $\kappa_c$ parameters where we have partially connected (i.e., $\bsG_{i,j}=0$ for $(i,j) \in\{(1,3),(3,1)\} $ and $\bsG_{i,j}=1$, otherwise) and fully connected (i.e., $\bsG_{i,j}=1$ for all $(i,j)$) populations. We can also see the average production levels $(\bar{X}^i_t)_{t\in[0,T]}$ for each population $i\in\{1,2,3\}$ in Fig.~\ref{fig:major_connection_effect}. We see that when the major player cares about the unhappiness of the people in their population more (i.e., when $\kappa_c$ is higher), they set a higher equilibrium carbon tax level. This can be explained by the fact that, when there is major player, producers care about their self interests and they may not care about the adverse effects of their carbon emission; however, major players know that higher carbon emission levels are worse for the society as a whole and they set carbon tax levels accordingly.

        \begin{table}
        \centering
        \caption{\small M4FNE Carbon tax levels for each population under different $\kappa_c$ parameters and different population connection setups.}
        \begin{tabular}{l |c c c |c c c}
        & \multicolumn{3}{c}{Partially Connected} & \multicolumn{3}{c}{Fully Connected}\\
        \hline \\[-2ex]
        $\kappa_c-1$ & $\hat{\gamma}^1$ & $\hat{\gamma}^2$ & $\hat{\gamma}^3$ & $\hat{\gamma}^1$ & $\hat{\gamma}^2$ & $\hat{\gamma}^3$ \\
        \hline
        $0.001$ & 2.6 & 3.7 & 1.1  & 5.4 & 3.8 & 2.2 \\
        $ 0.005$  & 12.1 & 17.4 & 4.5  & 26.1 & 17.8 & 9.5 \\
        $0.01$  & 24.1 & 34.6 & 8.6  & 52.1 & 35.3 & 18.4 \\
        \end{tabular}
        \label{tab:MajMin-NE_tax}%
        \end{table}%

\section{Conclusions and Future Work}
\label{sec:conclusion-future}

\subsection{Conclusions}

In this paper, we first introduce a general finite player game with multiple populations involving minor players and major players and give the Nash equilibrium definition both only for minor players in the populations and also for the full model with major players. Then, we discuss the multi-population major minor mean field game formulation of this problem in order to be able to have tractable solutions when number of players in the populations goes to infinity. Later, we give our motivating example with large number of electricity producers in $M$ populations and major players and solve it explicitly. We further discuss our results with a numerical example.

\subsection{Future Work}

For our future work, we have three main directions: The first direction is related to modeling. First, we plan to add the populations of consumers and second, we plan to add network structures inside the populations. The second future direction is related to experiments: we are interested in high scale experiments where we have large number of populations with more complex connections. The last future direction is related to exploring different types of equilibria in the general setup. For example, we plan to work on the Stackelberg multi-population mean field game equilibrium between the principals (i.e., major players) and mean field populations of minor players.

\bibliographystyle{ieeetr}


\begin{thebibliography}{10}

\bibitem{lasry2007mean}
J.-M. Lasry and P.-L. Lions, ``Mean field games,'' {\em Japanese journal of
  mathematics}, vol.~2, no.~1, pp.~229--260, 2007.

\bibitem{huang2006large}
M.~Huang, R.~P. Malham{\'e}, P.~E. Caines, {\em et~al.}, ``{Large population
  stochastic dynamic games: closed-loop McKean-Vlasov systems and the Nash
  certainty equivalence principle},'' {\em Communications in Information \&
  Systems}, vol.~6, no.~3, pp.~221--252, 2006.

\bibitem{cardaliaguet2018mean}
P.~Cardaliaguet and C.-A. Lehalle, ``Mean field game of controls and an
  application to trade crowding,'' {\em Mathematics and Financial Economics},
  vol.~12, pp.~335--363, 2018.

\bibitem{bensoussan2013mean}
A.~Bensoussan, J.~Frehse, P.~Yam, {\em et~al.}, {\em Mean field games and mean
  field type control theory}, vol.~101.
\newblock Springer, 2013.

\bibitem{carmona2018probabilistic}
R.~Carmona and F.~Delarue, {\em Probabilistic Theory of Mean Field Games with
  Applications I: Mean Field FBSDEs, Control, and Games}.
\newblock Probability Theory and Stochastic Modelling, Springer International
  Publishing, 2018.

\bibitem{cirant2015multi}
M.~Cirant, ``Multi-population mean field games systems with neumann boundary
  conditions,'' {\em Journal de Math{\'e}matiques Pures et Appliqu{\'e}es},
  vol.~103, no.~5, pp.~1294--1315, 2015.

\bibitem{bensoussan2018meanmulti}
A.~Bensoussan, T.~Huang, and M.~Lauri{\`e}re, ``Mean field control and mean
  field game models with several populations,'' {\em Minimax Theory and its
  Applications}, vol.~3, no.~2, pp.~173--209, 2018.

\bibitem{nourian2013epsilon}
M.~Nourian and P.~E. Caines, ``$\epsilon$-{N}ash mean field game theory for
  nonlinear stochastic dynamical systems with major and minor agents,'' {\em
  SIAM Journal on Control and Optimization}, vol.~51, no.~4, pp.~3302--3331,
  2013.

\bibitem{carmonazhu2016probabilistic}
R.~Carmona and X.~Zhu, ``A probabilistic approach to mean field games with
  major and minor players,'' 2016.

\bibitem{advertisement}
R.~Carmona and G.~Dayan\i{}kl\i{}, ``Mean field game model for an advertising
  competition in a duopoly,'' {\em International Game Theory Review}, vol.~23,
  no.~04, p.~2150024, 2021.

\bibitem{elie2019tale}
R.~Elie, T.~Mastrolia, and D.~Possama{\"\i}, ``A tale of a principal and many,
  many agents,'' {\em Mathematics of Operations Research}, vol.~44, no.~2,
  pp.~440--467, 2019.

\bibitem{dayanikli2023machine}
G.~Dayan{\i}kl{\i} and M.~Lauri\`ere, ``A machine learning method for
  stackelberg mean field games,'' 2023.

\bibitem{lachapelle2016efficiency}
A.~Lachapelle, J.-M. Lasry, C.-A. Lehalle, and P.-L. Lions, ``Efficiency of the
  price formation process in presence of high frequency participants: a mean
  field game analysis,'' {\em Mathematics and Financial Economics}, vol.~10,
  pp.~223--262, 2016.

\bibitem{gomessaude2021meanprice}
D.~A. Gomes and J.~Sa{\'u}de, ``A mean-field game approach to price
  formation,'' {\em Dynamic Games and Applications}, vol.~11, no.~1,
  pp.~29--53, 2021.

\bibitem{carmonafouquesun2015mean}
R.~Carmona, J.-P. Fouque, and L.-H. Sun, ``Mean field games and systemic
  risk,'' {\em Communications in Mathematical Sciences}, vol.~13, no.~4,
  pp.~911--933, 2015.

\bibitem{elie2020contact}
R.~Elie, E.~Hubert, and G.~Turinici, ``Contact rate epidemic control of
  covid-19: an equilibrium view,'' {\em Mathematical Modelling of Natural
  Phenomena}, vol.~15, p.~35, 2020.

\bibitem{aurell2022optimal}
A.~Aurell, R.~Carmona, G.~Dayan{\i}kl{\i}, and M.~Lauri\`ere, ``Optimal
  incentives to mitigate epidemics: a stackelberg mean field game approach,''
  {\em SIAM Journal on Control and Optimization}, vol.~60, no.~2,
  pp.~S294--S322, 2022.

\bibitem{alasseur2020extended}
C.~Alasseur, I.~Ben~Taher, and A.~Matoussi, ``An extended mean field game for
  storage in smart grids,'' {\em Journal of Optimization Theory and
  Applications}, vol.~184, pp.~644--670, 2020.

\bibitem{gueant_2010}
O.~Gu{\'e}ant, J.-M. Lasry, and P.~L. Lions, ``{Mean Field Games and oil
  production},'' in {\em {The Economics of Sustainable Development}}
  (Economica, ed.), 2010.

\bibitem{sircar_2017}
P.~Chan and R.~Sircar, ``Fracking, renewables, and mean field games,'' {\em
  SIAM Review}, vol.~59, pp.~588--615, 01 2017.

\bibitem{malhame_2017}
O.~Bahn, A.~Haurie, and R.~Malhame, {\em Limit Game Models for Climate Change
  Negotiations}, pp.~27--47.
\newblock 07 2017.

\bibitem{carbon_StackelbergMFG}
R.~Carmona, G.~Dayan{\i}kl{\i}, and M.~Lauri{\`e}re, ``Mean field models to
  regulate carbon emissions in electricity production,'' {\em Dynamic Games and
  Applications}, 2022.

\bibitem{shrivats2021principal}
A.~Shrivats, D.~Firoozi, and S.~Jaimungal, ``Principal agent mean field games
  in rec markets,'' {\em arXiv preprint arXiv:2112.11963}, 2021.

\bibitem{escribe2023mean}
C.~Escribe, J.~Garnier, and E.~Gobet, ``A mean field game model for renewable
  investment under long-term uncertainty and risk aversion,'' 2023.

\bibitem{gomes2016extended}
D.~A. Gomes and V.~K. Voskanyan, ``Extended deterministic mean-field games,''
  {\em SIAM Journal on Control and Optimization}, vol.~54, no.~2,
  pp.~1030--1055, 2016.

\bibitem{kobeissi2022classical}
Z.~Kobeissi, ``On classical solutions to the mean field game system of
  controls,'' {\em Communications in Partial Differential Equations}, vol.~47,
  no.~3, pp.~453--488, 2022.

\end{thebibliography}

\clearpage
\appendix

\subsection{Proof of Theorem~\ref{the:minors-FBSDE}}

We start with the necessity condition. Assume that the control profile $\hat\bsalpha$ is the equilibrium in the multi-population mean field game given the major players control profile $\bsgamma = (\gamma^i)_i$. Now, without loss of generality assume that the representative player in population $i$ deviates from their best response and use $\hat{\alpha}^i +\epsilon \check{\alpha}^i$. We note that since the player is insignificant, the mean field does not deviate. Then, we have

\begin{equation}
    \begin{aligned}
    \label{eq:perturbedcost}
        \delta J[\hat{\alpha}^i, \check{\alpha}^i] &:= \Big[\dfrac{d}{d\epsilon} J(\hat{\alpha}^i +\epsilon \check{\alpha}^i;\bar{\boldsymbol{X}}, \bar{\bsalpha}, \gamma^i)\Big]_{\epsilon=0}\\[1mm]
        &\hskip-0.5cm=\EE\Big[\int_0^T2\gamma^i \hat{\alpha}_t^i \check{\alpha}_t^i dt \\
        &+\kappa^i \sum_{i^\prime=1}^M\bsG_{i, i^\prime} 2(\hat{X}_T^i-\bar{X}_T^{i^{\prime}})\check{X}_T^i
        -\delta^i\check{X}_T^i \Big]
    \end{aligned}
\end{equation}
where we have the dynamics $d\check{X}^i_t = \eta^i \check{\alpha}_t^i dt + \sigma^i d\check{W}^i_t$ with initial condition $\check{X}^i_0=0$. Now we introduce the adjoint variable with the following dynamics:
\begin{equation*}
\begin{aligned}
        &dY^i_t = Z^i_t d W^i_t, 
        \\[1mm]
        &Y^i_T = -\delta^i + 2\kappa^i (\bsG[i,:] \mathds{1}_M  X^i_T-\bsG[i,:]\bar{X}_T).
\end{aligned}
\end{equation*}
Plugging the adjoint variable in the perturbed cost function \eqref{eq:perturbedcost} and applying integration by parts, we have:
\begin{equation*}
    \begin{aligned}
        &\Big[\dfrac{d}{d\epsilon} J(\hat{\alpha}^i +\epsilon \check{\alpha}^i;\bar{\boldsymbol{X}}, \bar{\bsalpha}, \gamma^i)\Big]_{\epsilon=0}\\[1mm]
        &= \EE\Big[\int_0^T2\gamma^i \hat{\alpha}_t^i \check{\alpha}_t^i dt +Y_T^i\check{X}_T^i\Big]\\[1mm]
        &=\EE\Big[\int_0^T2\gamma^i \hat{\alpha}_t^i \check{\alpha}_t^i dt + \int_0^T Y_t^id\check{X}_t^i + \int_0^T \check{X}_t^i dY_t^i\Big]\\[1mm]
        &=\EE\Big[\int_0^T\big( 2\gamma^i \hat{\alpha}_t^i +Y_t^i \eta_t^i \big)\check{\alpha}_t^i dt\Big].
    \end{aligned}
\end{equation*}
Since we assumed $(\hat\alpha^i_t)_t$ is the best response for the representative player, the above equation should be equal to 0 for any given $\check{\alpha}_t^i$. Thus, we have the following optimality condition:
\begin{equation}
\label{eq:opt_cond_minor}
    2\gamma^i \hat{\alpha}_t^i +Y_t^i \eta_t^i = 0,
\end{equation}
which implies that the optimal control is given in terms of the adjoint variable by: $\hat{\alpha}_t^i =  - \frac{\eta^i}{2 \gamma^i} Y^i_t,\ \forall i \in\llbracket M \rrbracket$. 

For the sufficiency, we assume that $\bsalpha$ is the multi-population mean field game Nash equilibrium. Then, without loss of generality we focus on a representative player in population $i$. We want to show by using the FBSDE system in Theorem~\ref{the:minors-FBSDE}, for all admissible $(\tilde{\alpha}^i_t)_t$, we have 
\begin{equation*}
    J(\alpha^i;\bar{\boldsymbol{X}}, \bar{\bsalpha}, \gamma^i) \leq J(\tilde{\alpha}^i;\bar{\boldsymbol{X}}, \bar{\bsalpha}, \gamma^i).
\end{equation*}
We have
\begin{equation*}
    \begin{aligned}
        &J(\alpha^i; \bar{\boldsymbol{X}}, \bar{\bsalpha}, \gamma^i) - J(\tilde{\alpha}^i;\bar{\boldsymbol{X}}, \bar{\bsalpha}, \gamma^i)\\
        &\hskip5mm=\EE\Big[\int_0^T \gamma^i[(\alpha_t^i)^2 - (\tilde{\alpha}_t^i)^2]dt \\[1mm]
        &\hskip9mm+ \kappa^i \sum_{i^{\prime}=1}^M \bsG_{i, i^{\prime}} [(X_T^i)^2 - (\tilde{X}_T^i)^2 -2 \bar{X}_T^{i^{\prime}}(X_T^i - \check{X}_T^i)]\\
        &\hskip9mm-\delta^i(X_T^i -\tilde{X}_T^i)\Big]\\
        &\hskip5mm =\EE\Big[\int_0^T \gamma^i[(\alpha_t^i)^2 - (\tilde{\alpha}_t^i)^2]dt\\[1mm]
        &\hskip9mm+ (X_T^i-\tilde X_T^i)Y_T^i - (X_T^i-\tilde X_T^i)^2\kappa^i\sum_{i^{\prime}=1}^M \bsG_{i, i^{\prime}}\Big]\\
        &\hskip5mm \leq \EE\Big[\int_0^T \gamma^i[2\alpha_t^i(\alpha_t^i - \tilde{\alpha}_t^i)]dt\\[1mm]
        &\hskip9mm+ (X_T^i-\tilde X_T^i)Y_T^i - (X_T^i-\tilde X_T^i)^2\kappa^i\sum_{i^{\prime}=1}^M \bsG_{i, i^{\prime}}\Big],
    \end{aligned}
\end{equation*}
where we used the subgradient inequality for the square function in the last line. Now by using integration by parts, we conclude:
\begin{equation*}
    \begin{aligned}
        &J(\alpha^i;\bar{\boldsymbol{X}}, \bar{\bsalpha}, \gamma^i) - J(\tilde{\alpha}^i;\bar{\boldsymbol{X}}, \bar{\bsalpha}, \gamma^i)\\
        &\hskip3mm \leq \EE\Big[\int_0^T \gamma^i[2\alpha_t^i(\alpha_t^i - \tilde{\alpha}_t^i)]dt+ (X_0^i-\tilde X_0^i)Y_0^i \\[1mm]
        &\hskip5mm + \int_0^T (X_t^i-\tilde X_t^i)dY_t^i + \int_0^T Y_t^id(X_t^i-\tilde X_t^i)\\[1mm]
        &\hskip5mm- (X_T^i-\tilde X_T^i)^2\kappa^i\sum_{i^{\prime}=1}^M \bsG_{i, i^{\prime}}\Big]\\[1mm]
        &\hskip3mm \leq \EE\Big[\int_0^T (2\gamma^i\alpha_t^i +Y_t^i \eta^i)(\alpha_t^i - \tilde{\alpha}_t^i)dt+ (X_0^i-\tilde{X}_0^i) Y_0^i + \\[1mm]
        &\hskip5mm + \int_0^T (X_t^i-\tilde X_t^i)Z_t^idW_t^i- (X_T^i-\tilde X_T^i)^2\kappa^i\sum_{i^{\prime}=1}^M \bsG_{i, i^{\prime}}\Big] \\[1mm]
       &\hskip3mm \leq 0.
    \end{aligned}
\end{equation*}
The last inequality follows from the optimality condition~\eqref{eq:opt_cond_minor}
and from the fact that both ${X}_0^i\sim \mu^i_0$ and $\tilde{X}_0^i\sim \mu^i_0$.

\hfill$\square$

\subsection{Proof of Proposition~\ref{prop:minors-ode}}

Theorem~\ref{the:minors-FBSDE} states that a solution of the FBSDE \eqref{eq:FBSDE_minor} gives the equilibrium in the mean field game given the major players' controls, $\bsgamma=(\gamma^i)_i$. In order to solve the FBSDE \eqref{eq:FBSDE_minor}, as mentioned in the main text, we introduce an ansatz for $Y^i_t$ as an affine function of $X^i_t$ in the form $Y^i_t = A^i_t  + B^i_t (\bsG[i,:] \mathds{1}_M X^i_t - \bsG[i,:]\bar{X}_t)$ for all $i \in \llbracket M\rrbracket$, which is consistent with the fact that $Y^i_t$ represents the derivative of the value function of a representative player in population $i$ and, as usual in linear-quadratic models, we expect the value function to be quadratic. Then by plugging in the ansatz, we have
\begin{equation*}
    \begin{aligned}
        dY_t^i =& \dot{A}_t^i dt + \dot{B}_t^i \Big(\bsG[i,:] \mathds{1}_M X^i_t - \bsG[i,:]\bar{X}_t\Big)dt\\
        & + {B}_t^i \Big(\bsG[i,:] \mathds{1}_M d X^i_t - \bsG[i,:]d\bar{X}_t\Big)\\
        =& \dot{A}_t^i dt + \dot{B}_t^i \Big(\bsG[i,:] \mathds{1}_M X^i_t - \bsG[i,:]\bar{X}_t\Big) dt \\
        & + {B}_t^i \Big(\bsG[i,:] \mathds{1}_M (-\frac{{|\eta^i|^2}}{2\gamma^i}Y^i_tdt +\sigma^i dW^i_t)\\
        & - \sum_{j=1}^M\bsG[i,j]\bar{Y}_t^jdt\Big)
    \end{aligned}
\end{equation*}
Then we plug in $dX^i_t = -\frac{{|\eta^i|^2}}{2\gamma^i}Y^i_tdt +\sigma^i dW^i_t = -\frac{{|\eta^i|^2}}{2\gamma^i}\Big(A^i_t  + B^i_t (\bsG[i,:] \mathds{1}_M X^i_t - \bsG[i,:]\bar{X}_t)\Big)dt +\sigma^i dW^i_t$ and $d\bar{X}^j_t = -\frac{{|\eta^j|^2}}{2\gamma^j}\bar{Y}^j_tdt = -\frac{{|\eta^j|^2}}{2\gamma^j}\big(A_t^j + B^j_t \big[\bsG[j,:]\mathds{1}_M\bar{X}^j_t - \bsG[j,:]\bar{X}_t\big] \big)$ and match the terms with $dY_t^i = Z_t^i dW_t^i$ to end up with the backward ODEs in the FBODE system. Furthermore, by using the proposed ansatz, we have $Y^i_T = A^i_T  + B^i_T (\bsG[i,:] \mathds{1}_M X^i_T - \bsG[i,:]\bar{X}_T)=-\delta^i + 2\kappa^i (\bsG[i,:] \mathds{1}_M  X^i_T-\bsG[i,:]\bar{X}_T)$. Therefore, we find the terminal conditions for the backward ODEs as $A^i_T = -\delta^i$, $B_T^i = 2\kappa^i$. The forward ODE comes directly from the dynamics of $\bar{X}^i_t$. We would like to stress that the final system consists of coupled forward backward equations for all populations $i\in \llbracket M\rrbracket$.

\hfill$\square$

\subsection{Proof of Theorem~\ref{thm:existence-uniqueness}}

\label{app:proof_of_minor_existence_uniqueness}

Showing there exists a unique mean field equilibrium given the major players' controls, $\bsgamma = (\gamma^i)_i$ corresponds to showing that there exists a unique solution to the FBODE system~\eqref{eq:FBODE_A}--\eqref{eq:FBODE_BarX} with the boundary conditions~\eqref{eq:FBODE_boundarycond}. 

We first realize that the differential equation for $B_t^i$ is not coupled with any of $A_t^i$, $\bar{X}_t^i$ or $B_t^j$ for all $j\neq i$ and can be solved by itself. Then, the unique solution is given as
\begin{equation*}
    B_t^i = \dfrac{2\kappa^i \gamma^i}{\gamma^i + \kappa^i |\eta^i|^2 \bsG[i,:]\mathds{1}_M(T-t)}, \forall t\in[0,T],\ \forall i\in\llbracket M\rrbracket.
\end{equation*}
Given $(B_t^i)_{t\in[0,T], i\in\llbracket M\rrbracket}$, now we want to show that there exists a unique solution to the forward backward ODE system for $(A^i_t)_{t\in[0,T], i\in\llbracket M\rrbracket}$ \eqref{eq:FBODE_A} and $(\bar{X}^i_t)_{t\in[0,T], i\in\llbracket M\rrbracket}$ \eqref{eq:FBODE_BarX}. In order to do so, we make use of Banach Fixed Point Theorem.

In order to simplify the notation, given $(B_t^i)_{t\in[0,T], i\in\llbracket M\rrbracket}$ we write the coupled forward backward ODE system in matrix form:
\begin{equation*}
    \begin{aligned}
        \dot{A}_t &= B_t (\tilde{\bsG}-\bsG)\Xi A_t + B_t \bsG\Xi B_t (\tilde{\bsG}-\bsG)\bar{X}_t,\\
        \dot{\bar{X}}_t&= - \Xi(A_t + B_t (\tilde{\bsG}-\bsG)\bar{X}_t),
    \end{aligned}
\end{equation*}
where $A_t := [A_t^1, \dots, A_t^M]^{\top}$, $\bar{X}_t := [\bar{X}_t^1, \dots, \bar{X}_t^M]^{\top}$, $B_t$ is a diagonal matrix where we have $B_t^i$ on the $i^{\rm th}$ diagonal entry, $\Xi$ is a diagonal matrix where we have $\frac{|\eta^i|^2}{2\gamma^i}$ on the $i^{\rm th}$ diagonal entry and $\tilde{\bsG}$ is a diagonal matrix where we have $\sum_{j=1}^M\bsG[i, j] = \bsG[i, :]\mathds{1}_M$ on the $i^{\rm th}$ diagonal entry. Now we can follow similar steps to the proof of~\cite[Theorem 5]{carbon_StackelbergMFG}:

We first fix $(A_t^{\prime})_t$ and $(A_t^{\prime\prime})_t$ and find the corresponding $(\bar{X}_t^{\prime})_t$ and $(\bar{X}_t^{\prime\prime})_t$ by solving the following ODEs:
\begin{equation*}
    \begin{aligned}
        d\bar{X}_t^{\prime} &= \big[- \Xi A^{\prime}_t - \Xi  B_t (\tilde{\bsG}-\bsG)\bar{X}^{\prime}_t\big]dt,\\[1mm]
        d\bar{X}_t^{\prime\prime} &= \big[- \Xi A^{\prime\prime}_t - \Xi  B_t (\tilde{\bsG}-\bsG)\bar{X}^{\prime\prime}_t\big]dt,
    \end{aligned}
\end{equation*}
with $\bar{X}_0^{\prime} = \bar{X}_0^{\prime\prime} = \bar{x}_0$. Now, let $\tilde{A_t}:=A_t^{\prime}-A_t^{\prime\prime}$ and $\tilde{\bar{X}}_t:= \bar{X}_t^{\prime} - \bar{X}_t^{\prime\prime}$. Then,
\begin{equation*}
    d\tilde{\bar{X}}_t = \big[- \Xi \tilde{A}_t - \Xi  B_t (\tilde{\bsG}-\bsG)\tilde{\bar{X}}_t\big]dt,
\end{equation*}
with $\tilde{\bar{X}}_0 = 0$. Next, we introduce $\check{A}^{\prime}_t$ and $\check{A}^{\prime\prime}_t$ that solve
\begin{equation*}
    \begin{aligned}
        d\check{A}^{\prime}_t &= \big[B_t (\tilde{\bsG}-\bsG)\Xi A^{\prime}_t + B_t \bsG\Xi B_t (\tilde{\bsG}-\bsG)\bar{X}^{\prime}_t\big]dt,\\[1mm]
        d\check{A}^{\prime\prime}_t &= \big[B_t (\tilde{\bsG}-\bsG)\Xi A^{\prime\prime}_t + B_t \bsG\Xi B_t (\tilde{\bsG}-\bsG)\bar{X}^{\prime\prime}_t\big]dt,
    \end{aligned}
\end{equation*}
where $\check{A}^{\prime}_T = \check{A}^{\prime}_T = -[\delta^1, \dots, \delta^M]^{\top}=:-\boldsymbol{\delta}$. We introduce $\check{\tilde{A}}_t:=\check{A}_t^{\prime}-\check{A}_t^{\prime\prime}$ that solves the following ODE
\begin{equation*}
    d\check{\tilde{A}}_t = \big[B_t (\tilde{\bsG}-\bsG)\Xi \tilde{A}_t + B_t \bsG\Xi B_t (\tilde{\bsG}-\bsG)\tilde{\bar{X}}_t\big]dt,
\end{equation*}
with $\check{\tilde{A}}_T=\boldsymbol{0}$. Now, we would like to show the mapping $A\mapsto \check{A}$ is a contraction mapping and then use the Banach Fixed Point Theorem to show the existence and uniqueness of the solution.

\noindent \textbf{Step 1:} First, we start by writing the dynamics for $\Vert \tilde{\bar{X}}_t \Vert^2$ where $\Vert \cdot\Vert$ denotes the 2-norm:
\begin{equation*}
\begin{aligned}
    d\Vert \tilde{\bar{X}}_t\Vert^2&=2|(\tilde{\bar{X}}_t)^{\top}  d \tilde{\bar{X}}_t |\\
    &= 2\big|(\tilde{\bar{X}}_t)^{\top} \big[- \Xi \tilde{A}_t - \Xi  B_t (\tilde{\bsG}-\bsG)\tilde{\bar{X}}_t\big]\big|dt.
\end{aligned}
\end{equation*}
By using the dynamics above and Young's inequality, we bound $\Vert \tilde{\bar{X}}_t \Vert^2$ as follows:
\begin{equation*}
    \begin{aligned}
        \Vert \tilde{\bar{X}}_t\Vert^2 &= \int_0^t 2\big|(\tilde{\bar{X}}_s)^{\top} \big[- \Xi \tilde{A}_s - \Xi  B_s (\tilde{\bsG}-\bsG)\tilde{\bar{X}}_s\big]\big|ds \\
        &\leq \int_0^t 2\Vert \Xi\Vert |<\tilde{\bar{X}}_s, \tilde{A}_s>| ds \\
        &\hskip5mm +\int_0^t 2\Vert \Xi  B_s (\tilde{\bsG}-\bsG) \Vert \Vert \tilde{\bar{X}}_s \Vert^2 ds\\ 
        & \leq \int_0^t \Vert \Xi\Vert \left(\Vert\tilde{\bar{X}}_s\Vert^2+ \Vert\tilde{A}_s\Vert^2\right) ds \\
        &\hskip5mm +\int_0^t 2\Vert \Xi  B_s (\tilde{\bsG}-\bsG) \Vert \Vert \tilde{\bar{X}}_s \Vert^2 ds,
    \end{aligned}
\end{equation*}
where we use the matrix norm induced by 2-norm which is smaller than the vector 2-norm. By using Gr\"{o}nwall's inequality, we have
\begin{equation}
    \begin{aligned}
    \label{eq:norm_barX_bound} 
        \Vert \tilde{\bar{X}}_t\Vert^2 &\leq e^{\int_0^t\big( 2\Vert \Xi  B_s (\tilde{\bsG}-\bsG) \Vert \Vert + \Vert \Xi\Vert\big)ds} \int_0^t \Vert \Xi\Vert \Vert\tilde{A}_s\Vert^2ds \\
        &\leq C_1 \int_0^T \Vert\tilde{A}_s\Vert^2 ds
    \end{aligned}
\end{equation}
where 
\begin{equation}
\label{eq:C1}
    C_1 = e^{T\big( 2\Vert \Xi\Vert  \Vert B\Vert_T \Vert(\tilde{\bsG}-\bsG) \Vert + \Vert \Xi\Vert\big)}\Vert \Xi\Vert
\end{equation}
with $\Vert B \Vert_T = \sup_{0\leq s \leq T} \Vert B_s \Vert$.

\noindent \textbf{Step 2:} Now, similarly to \textbf{Step 1}, we write the dynamics for $\Vert \check{\tilde{A}}_t \Vert^2$:
\begin{equation*}
    \begin{aligned}
         &d\Vert \check{\tilde{A}}_t \Vert^2 
         \\
         &\hskip3mm= 2 (\check{\tilde{A}}_t)^{\top}\big[B_t (\tilde{\bsG}-\bsG)\Xi \tilde{A}_t + B_t \bsG\Xi B_t (\tilde{\bsG}-\bsG)\tilde{\bar{X}}_t\big]dt
    \end{aligned}
\end{equation*}
By using the dynamics above and Young's inequality, we bound $\Vert \check{\tilde{A}}_t \Vert^2$ as follows:
\begin{equation*}
    \begin{aligned}
        &\Vert \check{\tilde{A}}_t \Vert^2 \\
        &= \int_t^T 2 (\check{\tilde{A}}_s)^{\top}\big[B_s ({\bsG}-\tilde{\bsG})\Xi \tilde{A}_s + B_s \bsG\Xi B_s (\bsG-\tilde{\bsG})\tilde{\bar{X}}_s\big]ds\\
        &\leq \int_t^T  \Vert B_s ({\bsG}-\tilde{\bsG})\Xi \Vert 2<\check{\tilde{A}}_s,\tilde{A}_s>ds\\
        &\hskip5mm+ \int_t^T \Vert B_s \bsG\Xi B_s (\bsG-\tilde{\bsG})\Vert 2<\check{\tilde{A}}_s,\tilde{\bar{X}}_s>ds\\
        &\leq \int_t^T \big( \Vert B_s ({\bsG}-\tilde{\bsG})\Xi \Vert + \Vert B_s \bsG\Xi B_s (\bsG-\tilde{\bsG})\Vert \big) \Vert \check{\tilde{A}}_s \Vert^2ds\\
        &\hskip5mm+ \int_t^T  \Vert B_s ({\bsG}-\tilde{\bsG})\Xi \Vert \Vert \tilde{A}\Vert^2ds\\
        &\hskip5mm+ \int_t^T   \Vert B_s \bsG\Xi B_s (\bsG-\tilde{\bsG})\Vert   \Vert \tilde{\bar{X}}_s\Vert^2ds.
    \end{aligned}
\end{equation*}
Now, we plug in the bound~\eqref{eq:norm_barX_bound} found for $\Vert \tilde{\bar{X}}_s\Vert^2$:
\begin{equation*}
    \begin{aligned}
        &\Vert \check{\tilde{A}}_t \Vert^2 \\
        &\leq \int_t^T \big( \Vert B_s ({\bsG}-\tilde{\bsG})\Xi \Vert + \Vert B_s \bsG\Xi B_s (\bsG-\tilde{\bsG})\Vert \big) \Vert \check{\tilde{A}}_s \Vert^2ds\\
        &\hskip5mm+ \int_0^T  \Vert B_s ({\bsG}-\tilde{\bsG})\Xi \Vert \Vert \tilde{A}\Vert^2ds\\
        &\hskip5mm+ \int_0^T   \Vert B_s \bsG\Xi B_s (\bsG-\tilde{\bsG})\Vert   C_1 \big(\int_0^T \Vert\tilde{A}_\tau\Vert^2 d \tau\big)ds\\
        &\leq \int_t^T \big( \Vert B_s ({\bsG}-\tilde{\bsG})\Xi \Vert + \Vert B_s \bsG\Xi B_s (\bsG-\tilde{\bsG})\Vert \big) \Vert \check{\tilde{A}}_s \Vert^2ds\\
        &\hskip5mm+ \int_0^T  \Vert B_s ({\bsG}-\tilde{\bsG})\Xi \Vert \Vert \tilde{A}_s\Vert^2ds\\
        &\hskip5mm+ T\Big( \Vert B\Vert_T \Vert\bsG\Xi\Vert \Vert B\Vert_T \Vert(\bsG-\tilde{\bsG})\Vert   C_1 \big(\int_0^T \Vert\tilde{A}_s\Vert^2 ds\big)\Big)\\
        &\leq e^{T\big(\Vert B\Vert_T \Vert({\bsG}-\tilde{\bsG})\Xi \Vert + \Vert B\Vert_T \Vert\bsG\Xi \Vert \Vert B\Vert_T \Vert(\bsG-\tilde{\bsG})\Vert\big)}\times\\
        &\int_0^T( \Vert B_s ({\bsG}-\tilde{\bsG})\Xi \Vert + T\Vert B\Vert^2_T \Vert\bsG\Xi\Vert  \Vert \bsG-\tilde{\bsG}\Vert C_1)\Vert\tilde{A}_s\Vert^2 ds.
    \end{aligned}
\end{equation*}
Now, we define $\Vert \tilde{A} \Vert_T:= \sup_{0\leq t \leq T} \Vert \tilde{A}_t \Vert$ and $\Vert \check{\tilde{A}} \Vert_T:= \sup_{0\leq t \leq T} \Vert \check{\tilde{A}}_t \Vert$. Then, we have $\Vert \check{\tilde{A}} \Vert_T\leq C_2\Vert \tilde{A} \Vert_T$ where $C_2$ is given as:
\begin{equation*}
    \begin{aligned}
        C_2 =& Te^{T\big(\Vert B\Vert_T \Vert({\bsG}-\tilde{\bsG})\Xi \Vert + \Vert B\Vert_T \Vert\bsG\Xi \Vert \Vert B\Vert_T \Vert(\bsG-\tilde{\bsG})\Vert\big)}\\
        &\times \Big(\Vert B\Vert_T \Vert({\bsG}-\tilde{\bsG})\Xi \Vert + T\Vert B\Vert^2_T \Vert\bsG\Xi\Vert  \Vert(\bsG-\tilde{\bsG})\Vert C_1\Big),
    \end{aligned}
\end{equation*}
where $C_1$ is given in~\eqref{eq:C1}. When $T$ is small enough, we have $C_2<1$ which yields that the mapping $A \mapsto \check{A}$ is a contraction mapping. Then by Banach Fixed Point Theorem, for any small enough $T>0$, there exists a unique solution to the forward backward ODE system, which in turn implies that there exists a unique mean field Nash equilibrium between the minor players given the major players' control profile $\bsgamma$.

\hfill$\square$

\subsection{Proof of Theorem~\ref{thm:majors-controls-formula}}

The best response of major player $i$, $\hat{\gamma}^i$, is found by minimizing her cost function given other major players' controls, $\bsgamma^{-i}$, and representative minor players' controls, $\bsalpha$:

\begin{equation*}
  \begin{aligned}  
    \hat{\gamma}^i =& \argmin_{\gamma^i}\ J^{0}(\gamma^i; \bsgamma^{-i}, \bar{\boldsymbol{X}}, \bar{\bsalpha})\\
    =& \argmin_{\gamma^i}\ \kappa_c \Big(    \EE\big[\int_0^T \big( \gamma^i ({\alpha}^i_t)^2
        + 
        \sum_{i'=1}^{M} \textbf{G}_{i,i'} (\bar \alpha_t^{i'})^2 
        \big) dt
        \\
        &
        + 
        \kappa^i\sum_{i'=1}^{M} \textbf{G}_{i,i'}({X}^{i}_T-\bar{X}_T^{i'})^2 -\delta^i X^{i}_T\big]\Big)
    \\&- \int_0^T \gamma^i\EE[(\alpha_t^i)^2] dt
     +\sum_{i'=1}^{M} \textbf{G}_{i,i'} (\gamma^i-\gamma^{i'})^2
    + \kappa_g (\gamma^i)^2\\
    =& \argmin_{\gamma^i}\ (\kappa_c-1) \int_0^T \gamma^i\EE[(\alpha_t^i)^2] dt \\
    & +\sum_{i'=1}^{M} \textbf{G}_{i,i'} (\gamma^i-\gamma^{i'})^2
    + \kappa_g (\gamma^i)^2,
\end{aligned}
\end{equation*}
where in the last inequality we dropped the terms that do not affect the minimizer. Since $\kappa_g>0$ and the entries of $\bsG$ are nonnegative, the last expression is convex in $\gamma^i$. We find the unique minimizer by taking the derivative of the last expression with respect to ${\gamma}^i$ and thus obtain the expression given in~\eqref{eq:major-best-response}. Realize that this expression is coupled with the given controls of the other major players, i.e., we can write $\hat{\gamma}^i = f^i(\bsgamma^{-i})$ where $f^i(\bsgamma^{-i})$ is the right hand side in~\eqref{eq:major-best-response}. In the equilibrium, the best response will be given to the other players' best responses and we can write
\begin{equation*}
    \begin{bmatrix}
        \hat{\gamma}^1\\
        \hat{\gamma}^2\\
        \vdots \\
        \hat{\gamma}^M
    \end{bmatrix}=
    \begin{bmatrix}
        f^1(\hat{\bsgamma}^{-1})\\
        f^2(\hat{\bsgamma}^{-2})\\
        \vdots \\
        f^M(\hat{\bsgamma}^{-M})
    \end{bmatrix}
\end{equation*}
By solving the matrix equations for $\hat{\bsgamma} = [\hat{\gamma}^1,\dots, \hat{\gamma}^M]^{\top}$, we conclude that the equilibrium behavior of major players given the controls of the representative minor players, $\bsalpha$, are given as in~\eqref{eq:major-eq-given-minor}.

\hfill$\square$

\end{document}